%
%
%
%
%
%
%
%
%
%
%
\scrollmode
\magnification=\magstep1
\parskip=\smallskipamount

\def\demo#1:{\par\medskip\noindent\it{#1}. \rm}
\def\ni{\noindent}               
\def\ll{\leftline}
\def\cl{\centerline}

\def\begin{\ll{}\vskip 10mm\nopagenumbers}  
\def\pn{\footline={\hss\tenrm\folio\hss}}   
\def\ii#1{\itemitem{#1}}

%
%
\outer\def\beginsection#1\par{\bigskip
  \message{#1}\leftline{\bf\&#1}
  \nobreak\smallskip\vskip-\parskip\noindent}

%
%
\outer\def\proclaim#1:#2\par{\medbreak\vskip-\parskip
    \noindent{\bf#1.\enspace}{\sl#2}
  \ifdim\lastskip<\medskipamount \removelastskip\penalty55\medskip\fi}

\def\endpr{\hfill $\spadesuit$ \medskip}

%
%
%
%


%
%
%
%

\def\C{{\bf C}}

\def\N{{\bf N}}

\def\R{{\bf R}}

\def\Z{{\bf Z}}


%
%
%
%

%
%
%
\def\a{\alpha}

\def\c{\chi}

\def\x{\xi}


%
%
%
%
\def\bs{\backslash}              
\def\di{\partial}                
\def\hra{\hookrightarrow}

%
%

\def\c*{{\bf C}^*}

%
%
\def\dim{{\rm dim}\,}                    
\def\holo{holomorphic}                   
\def\nbd{neighborhood}                   
\def\ra{real-analytic}                   

\def\tr{totally real}                    
\def\pc{polynomially convex}             


\def\iff{if and only if}

\def\hvb{holomorphic vector bundle}

\def\phe{proper holomorphic embedding}

\def\hci{holomorphic complete intersection}
\def\dci{differentiable complete intersection}

\def\Aut{{\rm Aut}}                         

\def\wt{\widetilde}

\begin
\cl{\bf ON COMPLETE INTERSECTIONS}
\bigskip
\cl{Franc Forstneri\v c}
\bigskip\medskip\rm

%
%
%
%

\beginsection 1. The results.

In this paper we give examples of closed complex submanifolds
in complex euclidean spaces which are \dci s but not 
\hci s (theorems 1.1 and 1.2). 
We also prove a result on removing intersections
of holomorphic mappings from Stein manifolds with certain
complex subvarieties in euclidean spaces;
theorem 1.3 below extends a result of Forster and 
Ramspott from 1966 [FRa].

Recall that a closed complex submanifold $Y$ of codimension $d$ in
a complex manifold $X$ is a {\it \hci} if there exist $d$ \holo\
functions $f_1,\ldots,f_d\in {\cal O}(X)$ such that
$$
    Y=\{x\in X\colon f_1(x)=\ldots=f_d(x)=0\}            \eqno(1)
$$
and the differentials $df_j(x)$ ($1\le j\le d$) are
$\C$-linearly independent at each point $x\in Y$.
These differential induce a trivialization of the
complex normal bundle $N_Y=TX|_Y/TY$ of $Y$ in $X$.
There is a partial converse when $X$ is a Stein manifold:
{\it If the normal bundle $N_Y$ is trivial then $Y$ is a \hci\
in some open \nbd\ of $Y$ in $X$}
(since a \nbd\ of $Y$ in $X$ is biholomorphic to a \nbd\ of the
zero section in the normal bundle $N_Y$; see [GR, p.\ 256]). 
Similarly, a smooth real submanifold $Y$ of real codimension 
$d$ in a smooth manifold $X$ is a {\it \dci\/} if there exist $d$
smooth real functions on $X$ satisfying (1), with
independent differentials along $Y$. 
For results on complete intersections we refer the
reader to the papers [Sch] and [BK] and the references therein.

\proclaim 1.1 Theorem:
There exists a three dimensional closed complex submanifold
in $\C^5$ which is a \dci\ but not a \hci.

More precisely, given any compact orientable two dimensional surface
$M$ of genus $g\ge 2$ we construct a three dimensional
Stein manifold $Y$ which is homotopically equivalent to $M$
and whose tangent bundle $TY$ is trivial as a real vector bundle,
but is non-trivial as a complex vector bundle. We then show that
any \phe\ $Y\hra \C^5$ (or $Y\hra\C^7$) satisfies the conclusion
of theorem 1.1. In fact, we prove

\medskip\ni\bf Theorem 1.2. \sl
Let $Y$ be a Stein manifold of dimension $m$ whose tangent
bundle is trivial as a real vector bundle, but is non-trivial
as a complex vector bundle. Choose integers $m$ and $d$ 
such that either
\item{(a)} $d=2$ and $m\in \{2,3\}$, or
\item{(b)} $d=4$ and $2\le m\le 7$.

\ni Then the image of any \phe\ $Y\hra\C^{m+d}$
is a \dci\ but not a \hci\ in $\C^{m+d}$.
\medskip\rm

Multiplying our $Y^3\subset \C^5$ by $\C^k$ we obtain similar
examples in higher dimensions. Submanifolds of this type don't
exist in $\C^n$ for $n\le 3$, but we don't know the answer
for two dimensional submanifolds in $\C^4$. Recall that every
smooth \holo\ curve in $\C^n$ is a \holo\ complete intersection
[FRa], and so is every complex hypersurface in $\C^n$
(since all divisors on $\C^n$ are principal).

%
%
\demo Example 1: There exists a Stein manifold $X$ of dimension
four and a closed complex submanifold $Y\subset X$ of dimension
two such that $Y$ is a \dci\ but not a \hci\ in $X$. We can
choose $Y$ to have the homotopy type of the real two-sphere.
(See proposition 2.4 in section 2).
\endpr

\pn

In the remainder of this section we discuss the problem of
removing intersections of \holo\ maps from Stein manifolds
into $\C^d$ with certain analytic subvarieties
$\Sigma\subset \C^d$. Our main result, theorem 1.3,
contains as a special case the result on complete intersections 
due to Forster and Ramspott [Fra]. To motivate the discussion 
we first look at the complete intersections problem in the
more general context of complex spaces (with singularities). 
A closed complex subvariety $Y$ of a complex space $X$ is a 
\hci\ in $X$ if there exist $d=\dim X-\dim Y$ global sections 
of the analytic sheaf of ideals ${\cal J}_Y$ which generate this 
sheaf at each point of $X$. Consider the short exact sequence
$$ 
   0\to {\cal J}^2_Y \to
   {\cal J}_Y \to {\cal J}_Y/{\cal J}_Y^2 \to 0.
$$
When $Y$ is a local complete intersection of codimension
$d$, the quotient $N^*_Y= {\cal J}_Y/{\cal J}_Y^2$ is a
locally trivial analytic sheaf of rank $d$ with support on $Y$,
that is, a \hvb\ of rank $d$ over $Y$. The dual bundle $N_Y$
of $N^*_Y$ is by definition the normal bundle of $Y$ in $X$;
in the smooth case this coincides with the usual definition
of $N_Y$.

Suppose now that $X$ is Stein and $Y\subset X$ is a local complete
intersection of codimension $d$ in $X$ with normal bundle $N_Y$.
If $Y$ is a complete intersection then $N_Y$ is trivial
(since its dual bundle $N_Y^*={\cal J}_Y/{\cal J}^2_Y$
is generated by the images of the generators of ${\cal J}_Y$
and hence is trivial). The following partial converse
was obtained in 1966 by Forster and Ramspott [FRa]
by using the Oka-Grauert homotopy principle ([Gra], [Car]):

%
%
\medskip\sl Let $Y$ be a local complete intersection
of codimension $d$ with trivial normal bundle in a Stein space
$X$. Suppose that $U\subset X$ is an open set containing $Y$ and
the functions $f=(f_1,\ldots,f_d)\in {\cal O}(U)^d$ generate
${\cal J}_Y$ on $U$. If there is a continuous map $\wt{f}\colon
X\to \C^d$ such that $\wt{f}=f$ near $Y$ and $\wt{f}^{-1}(0)=Y$,
then $Y$ is a \hci\ in $X$. Such $\wt{f}$ always exists if
$\dim Y< \dim X/2$, or if $X$ is contractible 
and $\dim Y\le 2(\dim X -1)/3$.
\medskip\rm

Furthermore, M.\ Schneider proved that for any local complete
intersection $Y\subset X$ with trivial normal bundle
the sheaf ${\cal J}_Y$ admits $d+1$ generators
(theorem 2.5 in [Sch]).

Suppose now that $\Sigma$ is a closed complex subvariety of $\C^d$,
$f\colon X\to\C^d$ a \holo\ map from a Stein manifold $X$
and $Y\subset X$ a connected component (or a union of
such components) of $f^{-1}(\Sigma)$. 
{\it When is it possible to modify $f$ to a \holo\ map
$g\colon X\to \C^d$ such that $g^{-1}(\Sigma)=Y$
and $g-f$ vanishes to a given order on $Y$ ?}
A necessary condition is that we can modify $f$ to a
{\it continuous\/} map with the required properties, and
we are interested in the corresponding
{\it homotopy principle}. Example 2 below shows that we must
restrict the class of subvarieties $\Sigma$ to obtain positive
results. Denote by $\Aut\C^d$ the group of all \holo\
automorphisms of $\C^d$.

\proclaim Definition 1: A closed complex subvariety $\Sigma\subset
\C^d$ is said to be {\it tame\/} if there is a $\Phi\in \Aut\C^d$
such that
$  \Phi(\Sigma)\subset
   \Gamma = \{(z',z_d)\in \C^d \colon |z_d|\le 1+|z'|\}.
$

\rm  Every proper complex algebraic subvariety of $\C^d$
is tame. Conversely, a subvariety $\Sigma \subset \C^d$ of pure
dimension $d-1$ contained in $\Gamma$ is algebraic [Chi],
and hence $\Sigma^{d-1}\subset \C^d$ is tame \iff\ it
is equivalent to an algebraic subset by an automorphism
of $\C^d$. For discrete sets our notion of tameness coincides
with that of Rosay and Rudin [RR].

%
%
\medskip\ni \bf 1.3 Theorem. \rm (Removal of intersections.) \sl
Let $\Sigma$ be a closed complex analytic subvariety of $\C^d$
satisfying one of the following conditions:
\item{(a)} $\Sigma$ is tame and $\dim \Sigma\le d-2$;
\item{(b)} a complex Lie group acts holomorphically
and transitively on $\C^d\bs \Sigma$.

\noindent
Let $X$ be a Stein manifold, $f\colon X \to \C^d$ a \holo\ map and
$Y \subset X$ a union of connected components of $f^{-1}(\Sigma)$.
If there is a continuous map $\wt{f}\colon X\to \C^d$
which equals $f$ in a \nbd\ of $Y$ and satisfies
$\wt{f}^{-1}(\Sigma)=Y$, then for each $r\in \N$ there is
a \holo\ map $g\colon X\to \C^d$ such that $g^{-1}(\Sigma)=Y$
and $g-f$ vanishes to order $r$ along $Y$. Such $g$ 
always exists if $\dim X< 2(d-\dim \Sigma)$, 
or if $X$ is contractible and $\dim Y \le 2(d-\dim \Sigma-1)$.
\medskip \rm

The theorem on complete intersections [FRa] mentioned above
corresponds to the special case of theorem 1.3 with 
$\Sigma=\{0\} \subset \C^d$; in this case we have
$\dim X=\dim Y+d$, $\dim \Sigma=0$, and hence the 
dimension conditions in both theorems agree.

Theorem 1.3 is proved in section 4. Using the 
Oka-Grauert-Gromov homotopy principle from [Gro], [FP1], [FP2] 
we reduce it to an extension problem for continuous maps to which 
we then apply some standard results from the obstruction theory. 
The following example shows that theorem 1.3 fails for non-tame 
subvarieties of $\C^d$, independently of their codimension.

\demo Example 2: For each $d\ge 1$ there is a discrete
set $\Sigma\subset \C^d$ such that every \holo\ map
$g\colon \C^d\to \C^d\bs \Sigma$ has rank at most $d-1$.
When $d=1$, this holds already if $\Sigma$ contains
two points (the complement is then hyperbolic);
for $d>1$ such sets were constructed by Rosay and Rudin [RR]. 
For such $\Sigma$ the conclusion of theorem 1.3 fails for 
$Y=\{0\}\subset \C^d=X$.
\endpr

We observe that complements of tame subvarieties of
codimension at least two admit Fatou-Bieberbach domains;
for proof see section 4:

\proclaim  1.4 Proposition: For each tame complex subvariety
$\Sigma\subset \C^d$ of codimension at least two there
exists an injective \holo\ map $F\colon \C^d\to \C^d\bs \Sigma$
(a Fatou-Bieberbach map). If $0\notin \Sigma$, we can choose
$F$ such that $F(0)=0$ and $F$ is tangent to the identity
at $0$ to arbitrary finite order. The same is true if
$\Sigma$ is a compact subset of $\C^d$ whose polynomial
hull does not contain the origin.

\proclaim 1.5 Corollary: Let $\Sigma \subset \C^d\bs \{0\}$
be as in proposition 1.4. If $Y\subset X$ is a complete
intersection of codimension $d$ in a complex space $X$,
we can choose generators $f_1,\ldots,f_d$ of ${\cal J}_Y$
such that the map $f=(f_1,\ldots,f_d) \colon X\to \C^d$
avoids $\Sigma$.

\demo Proof: If $g=(g_1,\ldots,g_d)$ is any set of generators
for ${\cal J}_Y$ and $F$ satisfies proposition 1.4, the components
of the map $f=F\circ g \colon X\to \C^d$ are generators
of ${\cal J}_Y$ and we have $f(X)\subset \C^d\bs \Sigma$.
\endpr

We conclude this introduction by mentioning two
open problems.

\demo Problem 1 (Murthy): Let $Y\subset \C^n$
be a local \hci\ with trivial normal bundle. Is
$Y$ a complete intersection in $\C^n$ ?
In particular, is every closed complex submanifold
$Y\subset \C^n$ with trivial normal bundle a
\hci\ in $\C^n$ ? The first open case to consider
is five dimensional submanifolds in $\C^8$ [Sch].
The answer is negative for \dci s (example 1.1 in [BK]).

\demo Problem 2: If the answer to problem 1 is negative
in general, we may ask whether there exists a closed
complex submanifold $Y\subset \C^n$ with the following properties:
\item{(a)} the complex normal bundle of $Y$ in $\C^n$ is trivial,
\item{(b)} $Y$ is a \dci\ in $\C^n$, but
\item{(c)} $Y$ is not a \hci\ in $\C^n$.

The paper is organized as follows. In section 2 we collect
some preliminary material on vector bundles. In section 3
we prove theorems 1.1 and 1.2. In section 4 we prove
theorem 1.3 and proposition 1.4.

%
%
%
%
\beginsection 2. Preliminaries.

We begin by recalling some basic facts on real and complex
vector bundles over CW-complexes; the proofs can be found in [Hus].
The results concerning complex vector bundles remain true for
holomorphic vector bundles over Stein spaces in view of
the Oka--Grauert principle [Gra], [Car]  and the fact that
any $n$-dimensional Stein space is homotopy equivalent to an
$n$-dimensional CW-complex [Ham].

\def\Vect{{\rm Vect}}

We denote by $\Vect_\R^k(X)$ (resp.\ $\Vect_\C^k(X)$) the
topological isomorphism classes of real (respectively complex)
vector bundles of rank $k$ over a CW-complex $X$. If $X$ is a
Stein space then by Grauert's theorem [Gra], $\Vect_\C^k(X)$
coincides with the equivalence classes of \hvb s of rank $k$
over $X$. By ${\cal T}^k_\R$ (resp.\ ${\cal T}^k_\C$) we denote
the trivial real (respectively complex) vector bundle of rank
$k$ over a given base (which will always be clear from
the context).

\proclaim 2.1 Theorem: Let $X$ be an $n$-dimensional
CW-complex. The map
$$ \Vect_\R^k(X)\to \Vect_\R^{k+r}(X),\quad
   E\to E\oplus {\cal T}_\R^{r} \qquad (k,r\ge 1)
$$
is surjective if $k\ge n$ and is bijective if
$k\ge n+1$.

\proclaim 2.2 Theorem: Let $X$ be an $n$-dimensional
CW-complex. The map
$$ \Vect_\C^k(X)\to \Vect_\C^{k+r}(X),\quad
   E\to E\oplus {\cal T}_\C^{r} \qquad (k,r\ge 1)
$$
is surjective when $k\ge [n/2]$ and is bijective
when $k\ge [{n+1\over 2}]$. In particular, if $E\to X$
is a nontrivial complex vector bundle of rank
$k\ge [{n+1\over 2}]$, the bundle $E\oplus {\cal T}_\C^r$
is nontrivial for each $r\in \N$.

\demo Remark: Theorem 2.2 shows that any complete
intersection submanifold $Y$ in $\C^n$ is parallelizable,
since $TY\oplus N_Y=T\C^n|_Y={\cal T}_\C^n$ and
$N_Y$ trivial implies $TY$ trivial. Likewise,
any real submanifold $Y\subset \R^N$ which is a
\dci\ is stably parallelizable, i.e.,
$TY\oplus{\cal T}_\R^1$ is trivial.
\endpr

We shall also need the following result from [BK].

\proclaim 2.3 Theorem: Each smooth submanifold $Y\subset \R^n$
of codimension $d \in \{1,2,4,8\}$ and with trivial normal
bundle is a \dci\ in $\R^n$.

\demo Proof: 
We recall the proof from [BK] for the sake of completeness. 
By triviality of $N_Y$ there is an open set $U\subset \R^n$ 
containing $Y$ and a smooth map $f=(f_1,\ldots,f_d)\colon U\to \R^d$
which defines $Y$ as a smooth complete intersection in $U$.
Let $U^*=U\bs Y$ and let $\phi\colon U^*\to S^{d-1}$
(the unit sphere in $\R^d$) be defined by
$\phi(x)=f(x)/||f(x)||$. If $d \in  \{2,4,8\}$,
$S^{d-1}$ admits $d-1$ linearly independent vector fields
$v_2,\ldots,v_d$. For $x\in U^*$ we denote by $A(x)$
the $d\times d$ matrix whose first column is $f(x)$ and
the subsequent columns are $v_j\circ\phi(x)$, $2\le j\le d$.
Let $E\to \R^n$ be the smooth rank $d$ vector bundle
obtained by patching the trivial bundles over the open 
covering $(U,\R^n\bs Y)$ of $\R^n$ by the map 
$A\colon U\bs Y\to GL(d,\R)$. Since $f=Ae_1$ for 
$e_1=(1,0,\ldots,0)^t$, the maps $f$ and $e_1$ patch 
together to a global section $\wt{f}\colon \R^n\to E$ 
which has no zeros outside of $U$. Since every vector bundle 
over $\R^n$ is trivial, $\wt{f}$ gives rise to a smooth map 
$\R^n\to \R^d$ which defines $Y$ as a complete 
intersection in $\R^n$.
\endpr

\demo Remark: Theorem 2.3 holds (with the same proof) if we
replace $\R^n$ by any contractible smooth manifold. However, the
argument does not apply if $d\notin\{1,2,4,8\}$, and the authors
of [BK] conjectured that the conclusion holds only for
the indicated values of $d$.
\endpr

The following result justifies example 1 in the introduction.

\proclaim 2.4 Proposition:
There exists a Stein manifold $X$ of dimension
four and a closed complex submanifold $Y\subset X$ of dimension
two which is homotopy equivalent to the two-sphere such that
$Y$ is a \dci\ but not a \hci\ in $X$.

\demo Proof: We take $X$ to be the total space of a rank two
\hvb\ over a two dimensional Stein manifold $Y$ such that the
bundle is trivial as a real vector bundle but non-trivial
as a complex vector bundle over $Y$. Its zero section,
which we identify with $Y$, is then a \dci\ but not
a \hci\ in $X$. To obtain such a bundle we let $S$ be the Riemann
sphere and set $E=TS\oplus {\cal T}_\C^1 \to S$, where $TS$ is the
\holo\ tangent bundle of $S$. Since $TS$ is non-trivial
and the base has dimension two, theorem 2.2 shows that $E$ is
non-trivial as a complex vector bundle. However, as a real bundle we
have $E=(TS\oplus {\cal T}_\R^1)\oplus {\cal T}_\R^1$ which is trivial.
We now take $Y$ to be a Stein complexification of $S$, containing
$S$ as a maximal real submanifold, and we extend $E$ to
a holomorphic vector bundle $X\to Y$.
\endpr

It is instructive to carry out the above procedure
explicitly by defining a non-trivial complex structure
on the trivial rank four bundle ${\cal T}^4_\R \to S$.
The last part of the argument below is essentially the same
one which can be used to prove theorem 2.2.

\demo Explicit construction of a non-trivial complex structure
on ${\cal T}^4_\R$ over the $2$-sphere:
Let $x=(x_1,x_2,x_3,x_4)$ be
real coordinates on $\R^4$ and let $\{{\bf e}_j\colon 1\le j\le
4\}$ be the corresponding standard basis of $T_x\R^4$. Let
$S\subset \{0\}\times \R^3 \subset \R^4$ be the unit hypersurface
sphere in the hyperplane $x_1=0$, and let
$V=T\R^4|_S \cong S\times \R^4$.
We can equip $V$ with the structure of a rank 2
complex vector bundle over $S$ by choosing a map
$J\colon S\to GL(4,\R)$ satisfying $J_x^2=-Id$\ for each $x\in S$.
One such choice is $J^0_x {\bf e}_1={\bf e}_2$,
$J^0_x {\bf e}_3={\bf e}_4$; in this structure $V\to S$
is a trivial $\C$-vector bundle over $S$.
Another choice is obtained by starting with
$$ J_x {\bf e}_1=
   x_2{\bf e}_2 + x_3 {\bf e}_3 + x_4{\bf e}_4,
   \qquad x=(0,x_2,x_3,x_4)\in S.
$$
Let $V^1_x \subset V_x$ be the real 2-plane spanned by
${\bf e}_1$ and $J_x{\bf e}_1$, and let $V^2_x \subset V_x$
denote the orthogonal complement to $V^1_x$. Notice that
$V^2 =TS$ and hence it is nontrivial. Since $V^2$ is an
oriented plane bundle, we can choose an orientation
preserving $J_x \colon V^2_x \to V^2_x$, depending continuously
on $x\in S$ and such that $J_x^2=-Id$ on $V^2_x$.
(The choice is unique if we require that $J_x$ be orthogonal.)
We then extend $J_x$ by linearity to $V_x$.

We claim that the $\C$-bundle $(V,J)$ over $S$
is not equivalent to the trivial $\C$-bundle $(V,J^0)$.
Suppose on the contrary that there exists an equivalence
$A\colon S\to GL(4,\R)$ between the two bundles,
meaning that $AJ^0A^{-1}=J$. The group preserving
$J^0$ is precisely $GL(2,\C)$, and hence for any map
$B\colon S\to GL(2,\C)$ we  have
$$ J=AJ^0A^{-1}= ABJ^0B^{-1}A^{-1} = (AB)J^0 (AB)^{-1}. $$
We claim that we can choose $B$ such that
$B^{-1}A^{-1}{\bf e}_1={\bf e}_1$ on $S$.
Since $\R^4\bs \{0\} \simeq S^3$, every map
$S=S^2\to \R^4\bs \{0\}$ is homotopic to a constant
map. Thus there is a homotopy $v_t\colon S\to \R^4\bs \{0\}$
($t\in [0,1]$) from $v_0={\bf e}_1$ to
$v_1=A^{-1}{\bf e}_1$. Denote by
$\tau\colon GL(2,\C)\to \R^4\bs \{0\}$ the map
$\tau(B)=B{\bf e}_1$. Clearly this map is a Serre
fibration, i.e., it has the homotopy lifting property.
Thus there is a homotopy $B_t\colon S\to GL(2,\C)$
($t\in [0,1]$), with $B_0=Id$, satisfying
$B_t{\bf e}_1=v_t$ for each $t\in [0,1]$.
At $t=1$ we get the desired map $B=B_1$ satisfying
$B{\bf e}_1=A^{-1}{\bf e}_1$.

Write $C=AB$; hence $J=CJ^0 C^{-1}$. By construction
we have $C{\bf e}_1={\bf e}_1$ and
$C{\bf e}_2 =CJ^0 {\bf e}_1= J{\bf e}_1$.
Thus $C$ maps the trivial subbundle
$U= \R^2\times \{0\}^2 \subset V$ onto the subbundle
$V^1\subset V$, and hence it induces an isomorphism
of quotient bundles $V/U\cong V/V^1 \cong V^2 \cong TS$.
This is a contradiction since the first bundle is
trivial while the second is not.

%
%
%
%
\beginsection 3. Proof of theorems 1.1 and 1.2.

\demo Proof of theorem 1.2:
Let $F\colon Y\to \C^n$, $n=m+d$, be any \phe. We identify
$Y$ with the submanifold $F(Y)\subset \C^n$ and denote
by $N_Y$ its \holo\ normal bundle. By the Oka-Cartan theory we
have a \holo\ splitting $T\C^n|_Y \cong TY\oplus N_Y$ [GR].
Since $Y$ is a Stein manifold of dimension $m$, 
it is homotopy equivalent to a real $m$-dimenisonal
CW-complex. Since $TY$ is a trivial real bundle of rank $2m$
over $Y$, theorem 2.1 shows that its complement $N_Y$ is also 
real trivial provided that $2d>m$. Furthermore, the 
real codimension of $Y$ is $2d$ which is assumed to be 
either $4$ or $8$ and hence $Y$ is a \dci\ in $\C^n$ by 
theorem 2.3. On the other hand, since $TY$ is non-trivial as 
a complex bundle over $Y$, theorem 2.2 implies that
$N_Y$ is also a non-trivial complex vector bundle
and hence  $Y$ is not a \hci\ in any open set $U\subset \C^n$
containing $Y$. (There are no restrictions on $d$ and $m$ 
in the last argument).
\endpr

In the proof of theorem 1.1 we shall need the following:

\proclaim 3.1 Proposition: For any compact orientable two dimensional
surface $M$ of genus $g\ge 2$ there exists a three dimensional
Stein manifold which is homotopically equivalent to $M$
and whose tangent bundle is trivial as a real vector
bundle but is non-trivial as a complex vector bundle.

\demo Proof:
Let $M$ be any surface as in the proposition; such $M$
is the connected sum of $g\ge 2$ tori. Its tangent bundle 
$TM$ is non-trivial, but $TM\oplus {\cal T}_\R^1$ is 
trivial since $M$ embeds as a real hypersurface in $\R^3$.
By theorem 1.8 in [For1] there exists a smooth (even \ra)
embedding $M\hra\C^2$ which is \tr\ except at finitely
many complex tangent points which are hyperbolic in the 
sense of Bishop [Bis] and such that the embedded submanifold 
$M\subset \C^2$ has arbitrary small Stein \nbd s 
$\Omega \subset \C^2$ with a deformation retraction 
$\pi\colon \Omega \to M$. 

We endow $TM$ with the structure of a complex line bundle 
and take $E= \pi^*(TM)\to \Omega$. By the Oka--Grauert 
theorem [Gra] the bundle $p\colon E\to\Omega$ has an equivalent 
structure of a \hvb. In the present situation we can obtain such 
a structure quite explicitly as follows.
Assume (as we may) that the embedding $M\subset \C^2$
is \ra. We can represent the bundle $TM$ by a
1-cocycle defined by \ra\ functions
$c_{ij} \colon U_{ij}\to \C\bs \{0\}$ on a (finite)
open covering ${\cal U}=\{U_i\}$ of $M$ such that the
closure of each of the sets $U_{ij}=U_i\cap U_j$
for $i\ne j$ is contained in the \tr\ part of $M$
(we only need to avoid the finitely many complex tangents in $M$).
The complexifications of the functions $c_{ij}$ now determine
a \holo\ line bundle structure on $E$ over an open \nbd\ 
of $M$ in $\C^2$. 

We claim that the total space $E$ of this \hvb\ satisfies 
proposition 3.1. Since the base $\Omega$ is Stein, $E$ 
is also Stein. Clearly $E$ is homotopy equivalent to 
$\Omega$ and hence to $M$. We identify $\Omega$ with the 
zero section of $E$. The tangent bundle of $E$ equals 
$TE = p^*(TE|_\Omega)$, where
$TE|_\Omega =T\Omega\oplus E ={\cal T}_\C^2 \oplus E$.
Since $E\to \Omega$ is a non-trivial bundle and the base 
$\Omega$  is homotopic to the surface $M$, theorem 2.2 
shows that $TE|_\Omega$ is non-trivial as a complex vector 
bundle over $\Omega$, and hence $TE$ is a non-trivial complex 
vector bundle. On the other hand, as real vector bundles
we have $TE|_\Omega = {\cal T}_\R^4 \oplus E= 
{\cal T}_\R^3\oplus ({\cal T}_\R^1\oplus E)$.
We have already observed that the second summand is trivial
and hence $TE|_\Omega$ is a trivial real bundle. Therefore 
$TE$ is also trivial as a real bundle over $E$.
\endpr

\demo Proof of theorem 1.1: Let $Y$ be a Stein manifold 
satisfying proposition 3.1. By the embedding theorem of 
Eliashberg and Gromov [EGr] and Sch\"urmann [Sch\"ur] there 
exists a \phe\ $Y\hra\C^5$. By theorem 1.2 $Y$ is then 
a \dci\ but not a \hci\ in $\C^5$. The same argument applies 
to any embedding $Y\hra \C^7$.
\endpr

%
%
%
%
\beginsection 4. Removal of intersections.

\demo Proof of proposition 1.4:
Consider first the case when $\Sigma\subset \C^d$ is a tame
subvariety of dimension at most $d-2$.
For $1\le j\le d$ we denote by $\pi_j\colon \C^d\to \C^{d-1}_j$
the projection onto the coordinate hyperplane $\{z_j=0\}$.
Tameness of $\Sigma$ implies that, after a biholomorphic change
of coordinates on $\C^d$, the restriction of $\pi_j$ to $\Sigma$
is proper for each $j$, and hence $\Sigma_j=\pi_j(\Sigma)$
is a proper closed analytic subset of $\C^{d-1}_j$.
By translation we may assume that $\Sigma_j$ does not
contain the origin for any $j$. Choose a \holo\ function
$g_j$ on $\C^{d-1}_j$ such that $g_j(0)=-\log 2$ and $g_j=0$
on $\Sigma_j$, and set
$$ \Phi_j(z)=\left(
   z_1,\ldots,z_{j-1},e^{g_j(\hat z_j)}z_j,
   z_{j+1}, \ldots,z_d \right),
$$
where $\hat z_j=(z_1,\ldots,z_{j-1},z_{j+1},\ldots,z_n)$. Then
$\Phi=\Phi_1\circ\Phi_2\circ\cdots \circ \Phi_d \in \Aut\C^d$,
$\Phi$ restricts to  the identity on $\Sigma$, $\Phi(0)=0$,
and $D\Phi(0)={1\over 2}I$. Thus $0\in \C^d$ is an attracting
fixed point of $\Phi$ whose basin of attraction is a Fatou-Bieberbach
domain $\Omega\subset \C^d \bs \Sigma$. We obtain the corresponding
Fatou-Bieberbach map $F\colon \C^d\to \Omega$ as in [RR].

If $K$ is a compact subset in $\C^d$ whose polynomial
hull does not contain the origin, we can construct a
Fatou-Bieberbach map $F\colon \C^d\to \C^d\bs K$
by the push-out method of Dixon and Esterle [DE]
(see also [For2]). Here is the outline. Replacing $K$
by $\widehat K$ we may assume that $K$ is \pc\ and
$0\notin K$. Denote by $B_r$ the closed ball of
radius $r$ in $\C^d$. Proposition 2.1 in [For2]
(or the results in [FRo]) gives an automorphism
$G_0\in \Aut\C^d$ which is tangent to identity
to a given order $r$ at $0$ and satisfies
$G_0(K)\cap B_1=\emptyset$. Set $K_1=G_0(K)$. Next we choose
$\Phi_1\in \Aut\C^d$ which approximates the identity
map on $B_1$, is tangent to the identity at $0$ and
satisfies $\Phi_1(K_1)\cap B_2=\emptyset$, and we let
$G_1=\Phi_1\circ G_0$. Continuing inductively we obtain
a sequence $G_j\in \Aut\C^d$ ($j\in \Z_+$) which converges
on some domain $\Omega\subset \C^d$ to a biholomorphic
map $G\colon \Omega\to \C^d$ of $G$ {\it onto} $\C^d$
(proposition 5.1 in [For2]). By construction $G$
is tangent to the identity to order $r$ at $0$ and
$K\cap \Omega=\emptyset$. The map $F=G^{-1}\colon \C^d\to \Omega$
satisfies proposition 1.4.
\endpr

\demo Proof of theorem 1.3:
We use the same notation as in the statement of the theorem.
Let ${\cal J}_\Sigma$ be the sheaf of ideals of $\Sigma\subset \C^d$
and let ${\cal J}_Y$ be the sheaf of ideals of $Y \subset X$.
We define an analytic sheaf of ideals ${\cal S}$ on $X$
as follows: at points $x\in Y$ we take ${\cal S}_x$ to be
the pull-back of ${\cal J}_{\Sigma,f(x)}$ by $f$, and
for $x\in X\bs Y$ we take ${\cal S}_x={\cal O}_{X,x}$.
More precisely, if $x\in Y$ and if ${\cal J}_\Sigma$
is generated by functions $h_1,\ldots,h_m$ in some
\nbd\ of $f(x)$ in $\C^d$, we take the functions
$h_j\circ f$ ($1\le j\le m$) as the generators
of ${\cal S}$ in a \nbd\ of $x$. Clearly ${\cal S}$
is a coherent analytic sheaf of ideals on $X$ and
${\cal O}_X/{\cal S}$ is supported on $Y$.

Choose $r\in \N$ and let ${\cal R}={\cal S} {\cal J}_Y^r$;
this is also a coherent sheaf of ideals on $X$ which coincides
with ${\cal O}_X$ on $X\bs Y$. By the Oka-Cartan theory 
there are finitely many global sections 
$\x_1,\ldots,\xi_k$ of ${\cal R}$
such that $Y=\{x\in X\colon \xi_j(x)=0,\ 1\le j\le k\}$. (We do
not require that the $\xi_j$'s generate ${\cal R}$ !)
We seek a map $g\colon X\to \C^d$ satisfying theorem 1.3
in the form
$$ g(x)= f(x) + \sum_{j=1}^k \xi_j(x) g_j(x)
        = f(x) + G(x)\xi(x),                           \eqno(2)
$$
where $G(x)=(g_1(x),\ldots,g_k(x))$ is a \holo\ $d\times k$
matrix-valued function and $\xi=(\xi_1,\ldots,\xi_k)^t$.
For any choice of $G$ the map $g=f+G\xi$ agrees with $f$ to
order $r+1$ along $Y$. Our goal is to choose $G$ such that
$g^{-1}(\Sigma)=Y$. Define a \holo\ map
$\Phi\colon X\times \C^{dk}\to \C^d$ by
$$ \Phi(x,v_1,\ldots,v_k) =
   f(x) + \sum_{j=1}^k \xi_j(x) v_j  
   \qquad (x\in X,\ v_1,\ldots,v_k\in \C^d 
$$
and let
$$ \wt{\Sigma}=\Phi^{-1}(\Sigma) \bs (Y\times \C^{dk}). \eqno(3) $$
Then the map $g=f+G\xi$ satisfies theorem 1.3 \iff\ $G$ is 
holomorphic and its graph in $X\times \C^{dk}$ avoids 
$\wt{\Sigma}$ defined by (3).

Observe that for each fixed $x\in X\bs Y$ the map
$\Phi(x,\cdotp)\colon \C^{dk}\to\C^d$ is an affine
surjection, while for $x\in Y$ we have $\Phi(x,\cdotp)=f(x)$
(hence $Y\times \C^{dk} \subset \Phi^{-1}(\Sigma)$).
Let $p\colon X\times \C^{dk}\to X$ denote the base projection.
We shall need the following	lemma.

\proclaim 4.1 Lemma: The set $\wt{\Sigma}$ defined by (3) 
is a closed complex subvariety of $X\times \C^{dk}$. 
Moreover, for each point $a\in X\bs Y$ there is a \nbd\ 
$U\subset X\bs Y$ of $a$ and a biholomorphic self-map
$\Psi$ of $\wt{U}= U\times \C^{dk}$, with  $p\circ \Psi=p$,
such that $\Psi(x,\cdotp)$ is affine linear for each $x\in U$
and $\Psi(\wt{\Sigma} \cap \wt{U})=U\times (\Sigma\times\C^{d(k-1)})$.

\demo Proof:
By definition $\wt{\Sigma}$ is a closed complex
subvariety in $(X\bs Y)\times \C^{dk}$. The second
statement follows immediatelly from the observation that
$\Phi(x,\cdotp) \colon \C^{dk} \to \C^d$ is an affine
surjection for any $x\in X\bs Y$ and hence is locally
(with respect to the base) equivalent to the projection
of $\C^{dk}$ onto $\C^d\times \{0\}^{d(k-1)}$.

It remains to show that $\wt{\Sigma}$ is closed in
$X\times \C^{dk}$. We need to show that, as $x\in X\bs Y$
approaches a point $x_0\in Y$, the fibers $\wt{\Sigma}_x$
leave any compact subset of $\C^{dk}$.
Choose a \nbd\ $V\subset \C^d$ of the point $f(x_0)$
and \holo\ functions $h=(h_1,\ldots,h_m)$ on $V$
which generate the ideal sheaf ${\cal J}_\Sigma$ on $V$.
Also choose a \nbd\ $U\subset X$ of $x_0$ with
$f(U)\subset V$. Let $\xi_1,\ldots,\xi_k$ be sections
of the sheaf ${\cal R}$ as above. By Taylor expansion
of $h$ at the point $f(x)$ for $x\in U$ and 
$v=(v_1,\ldots,v_k) \in \C^{dk}$ we get
$$ \eqalign{ h\circ\Phi(x,v) &=
             h\left( f(x)+ \sum_{j=1}^k \xi_j(x)v_j \right) \cr
             &= h(f(x))+ \sum_{|\a|\ge 1}
                c_\a D^\a h(f(x))
                \left( \sum_{j=1}^k  \xi_j(x) v_j \right)^\a \cr
             &= h(f(x)) + A(x,v)\xi(x), \cr}
$$
where $A$ is a \holo\ $d\times k$ matrix function.
Denoting by $||\cdotp||$ the Euclidean norm on $\C^d$
(and the corresponding matrix norm) we have
$$ ||h(\Phi(x,v))|| \ge
   ||h(f(x))|| - ||A(x,v)||\, \cdotp ||\xi(x)||.
$$
The components of $h(f(x))$ generate the sheaf ${\cal S}$
at each point of $U$. Hence, as $x\to x_0\in Y$, the term
$||\xi(x)||$ is of size $o(||h(f(x))||)$ by the definition
of the sheaf ${\cal R}={\cal S}{\cal J}^r$.
Hence for each $C>0$ there is a \nbd\ $U_C \subset U$ of
$x_0$ such that for all $x\in U_C$ and $v\in \C^{dk}$ with
$||v||\le C$ we have
$||h\circ \Phi(x,v)|| \ge ||h(f(x))||/2$, and hence
$\Phi(x,v)\in \Sigma$ \iff\ $x\in Y$. Thus for $x\in U_C$
the fiber $\wt{\Sigma}_x$ does not intersect the ball
of radius $C$ in $\C^{dk}$. This proves that $\wt{\Sigma}$
is closed in $X\times \C^{dk}$.
\endpr

We continue with the proof of theorem 1.3. The assumptions
on $\Sigma$ imply that the complement $\C^d\bs \Sigma$ admits
a spray in the sense of Gromov (see [FP1] and lemma 7.1 in [FP2]).
From this and the second statement in lemma 4.1 it follows that
the \holo\ submersion $h\colon Z=(X\times \C^{dk})\bs \wt{\Sigma} \to X$
admits a fiber dominating spray in a small \nbd\ of
any point $x\in X\bs Y$ ([Gro] or definition 1.1 in [FP2]).
By theorem 1.2 in [FP2] (see also [Gro], 4.5 Main Theorem)
the homotopy principle holds for sections of $Z$, meaning
that any continuous section $\wt{G} \colon X\to Z$
can be deformed to a \holo\ section.

A continuous extension $\wt{f} \colon X\to \C^d$
of $f$ as in theorem 1.3 can be lifted to a continuous
section $\wt{G} \colon X\to Z$ which is \holo\ near $Y$
(see lemma 8.1 in [FP2]). The homotopy principle gives
a \holo\ section $G\colon X\to Z$ such that the corresponding
map $g \colon X\to \C^d$ (2) satisfies theorem 1.3.

In the remainder we investigate the existence of a continuous
extension  $\wt{f}$ using the obstruction theory (see e.g.\
section V.5 in [Whi]). By [Ham] the subvariety $Y$ has a closed
\nbd\ $A \subset X$ such that the pair $(X,A)$ is homotopy
equivalent to a relative CW-complex of dimension $n=\dim X$
and $Y$ is a deformation retraction of $A$. Moreover, we may
choose $A$ so small that $\{x\in A\colon f(x)\in \Sigma\}=Y$.
Hence $f$ maps $A^*=A\bs Y$ to $\Omega=\C^d\bs \Sigma$, and
we wish to find an extension of $f$ to a
map from $X^*=X\bs Y$ to $\Omega$.

The pair $(X^*,A^*)$ can be represented by the same relative
CW-complex as $(X,A)$. Denote by $X_q$ its $q$-dimensional
skeleton, so our goal is to extend $f$ to a map $X_n\to \Omega$.
We begin by extending $f$ to the zero-skeleton $X_0$ by
arbitrarily prescribing the values at the points of $X_0$.
Since $\Omega$ is connected, we can further extend to a
map $f_1\colon X_1 \to \Omega$. Suppose inductively that
$f$ has already been extended to $f_q\colon X_q\to \Omega$
for some $q\ge 1$. The next skeleton $X_{q+1}$ is
obtained by attaching $(q+1)$-cells $e_{q+1}$ to
$X_q$ by maps $\di e_{q+1}\to X_q$. Composing this attaching
map with $f_q\colon X_q\to \Omega$ we obtain for each such
cell $e_{q+1}$ a map $\di e_{q+1}\to \Omega$ which defines an
element of the fundamental group $\pi_q(\Omega)$. In this way
we obtain a singular cochain
$c^{q+1}\in \Gamma^{q+1}(X^*,A^*;\pi_q(\Omega))$
(which is in fact a $(q+1)$-cocycle, called the
{\it obstruction cocycle}), and $f_q$ extends
to a map $f_{q+1}\colon X_{q+1}\to \Omega$ \iff\
$c^{q+1}=0$.

In our case we have $\pi_q(\Omega)=\pi_q(\C^d\bs \Sigma)=0$
for $1\le q\le 2s-2$, where $s=d-\dim \Sigma$.
This implies that $f$ can be extended to the skeleton
$X_{2s-1}$. Hence, if $\dim X< 2s=2(d-\dim \Sigma)$,
we have an extension $\wt{f}\colon X\bs Y\to \C^d\bs \Sigma$
as required.

Assume now that $X$ is contractible (e.g., $X=\C^n$).
We shall use the following more precise result
from obstruction theory ([Whi], theorem V.5.14):

\medskip\ni\sl Let $f_q\colon X_q\to \Omega$ for some $q\ge 1$.
Then $f_q|X_{q-1}$ can be extended to a map
$f_{q+1} \colon X_{q+1}\to\Omega$ \iff\
$\gamma^{q+1}(f)=[c^{q+1}]=0\in H^{q+1}(X^*,A^*;\pi_q(\Omega))$,
i.e., the cohomology class of the obstruction cocycle $c^{q+1}$
equals zero.
\medskip\rm

By excision we have $H^q(X^*,A^*;G)=H^q(X,A;G)$ for any abelian
coefficient group $G$. Since $X$ is contractible, the long exact
sequence for the cohomology of the pair $A\hra X$ gives
$H^{q+1}(X,A;G)=H^q(A;G)$ for $q\ge 1$. Furthermore, since $Y$
is a deformation retract of $A$ we have $H^q(A;G)=H^q(Y;G)$.
Together we obtain
$$ H^{q+1}(X^*,A^*;\pi_q(\Omega))= H^q(Y; \pi_q(\Omega)),
 \quad q\ge 1.
$$
Since $Y$ is a Stein manifold of dimension $m$, it is homotopy
equivalent to an $m$-dimensional CW-complex and hence
$H^q(Y; \pi_q(\Omega))=0$ for $q>m$. Thus, if
$f \colon A^*\to\Omega$ admits an extension to the
skeleton $X_{m+1}$, it also admits an extension to all
higher dimensional skeleta and hence to $X^*$.
Earlier we have seen that there is an extension to
$X_{2s-1}$ with $s=d-\dim \Sigma$. If we assume
$m+1\le 2s-1$, we thus obtain a desired continuous
extension of $f$ to $X^*$. This completes the proof
of theorem 1.3.
\endpr

\demo Acknowledgements: I wish to thank E.\ L.\ Stout
for having brought to my attention the question answered by
theorem 1.1 and for stimulating discussions.
This research was supported in part by an NSF grant,
by the Vilas foundation at the University of Wisconsin--Madison,
and by the Ministry of Science of the Republic of Slovenia.

\vfill\eject
%
%
%
%
\medskip\ni\bf References. \rm

\ii{[Bis]} E.\ Bishop:
Differentiable manifolds in complex Euclidean space.
Duke Math.\ J.\ {\bf 32} (1965), 1--21.

\ii{[BK]} J.\ Bochnak and W.\ Kucharz:
Complete intersections in differential topology and analytic geometry.
Boll.\ Un.\ Mat.\ Ital.\ B (7) {\bf 10} (1996), no.\ 4, 1019--1041.

\ii{[Car]} H.\ Cartan: Espaces fibr\'es analytiques.
Symposium Internat.\ de topologia algebraica, Mexico, 97--121 (1958).
(Also in Oeuvres, vol.\ 2, Springer, New York, 1979.)

\ii{[Chi]} E.\ M.\ Chirka: Complex analytic sets.
Kluwer Academic Publishers Group, Dordrecht, 1989.

\ii{[DE]} P.\ G.\ Dixon and J.\ Esterle:
Michael's problem and the Poincar\'e--Fatou--Bieber\-bach phenomenon.
Bull.\ Amer.\ Math.\ Soc. (N.S.) {\bf 15} (1986), 127--187.

\ii{[EGr]} Y.\ Eliashberg, M.\ Gromov: Embeddings of Stein manifolds.
Ann.\ Math.\ {\bf 136}, 123--135 (1992).

\ii{[FRa]}  O.\ Forster and K.\ J.\ Ramspott:
Analytische Modulgarben und Endromisb\"undel.
Invent.\ Math.\ {\bf 2} (1966), 145--170.

\ii{[For1]} F.\ Forstneri\v c:
Complex tangents of real surfaces in complex surfaces.
Duke Math.\ J.\ {\bf 67} (1992), 353--376.

\ii{[For2]} F.\ Forstneri\v c:
Interpolation by holomorphic automorphisms and embeddings in $\C^n$.
J.\ Geom.\ Anal. {\bf 9}, no.1, (1999) 93-118.

\ii{[FP1]} F.\ Forstneri\v c and J.\ Prezelj:
Oka's principle for holomorphic fiber bundles with sprays.
Math.\ Ann.\ {\bf 317} (2000), 117-154.

\ii{[FP2]} F.\ Forstneri\v c and J.\ Prezelj:
Oka's principle for holomorphic submersions with sprays.
Preprint, 1999.

\ii{[FRo]} F.\ Forstneri\v c and J.-P.\ Rosay:
Approximation of biholomorphic mappings by automorphisms
of $\C^n$.
Invent.\ Math.\ {\bf 112} (1993), 323--349.
Erratum, \it Invent.\ Math.\ \bf 118 \rm (1994), 573--574.

\ii{[Gra]} H.\ Grauert: Analytische Faserungen \"uber
holomorph-vollst\"andigen R\"aumen. \break
Math.\ Ann.\ {\bf 135}  (1958), 263--273.

\ii{[GRe]} H.\ Grauert and R.\ Remmert: Theory of Stein Spaces.
Grundl.\ Math.\ Wiss.\ {\bf 227}, Springer, New York, 1977.

\ii{[Gro]} M.\ Gromov:
Oka's principle for holomorphic sections of elliptic bundles.
J.\ Amer.\ Math.\ Soc.\ {\bf 2}  (1989), 851-897.

\ii{[GR]} C.\ Gunning and H.\ Rossi:
Analytic functions of several complex variables.
Pren\-tice--Hall, Englewood Cliffs, 1965.

\ii{[Ham]} H.\ A.\ Hamm: Zur Homotopietyp Steinscher R\"aume.
J.\ Reine Angew.\ Math.\ {\bf 338} (1983), 121--135.

\ii{[Hus]} D.\ Husemoller: Fibre bundles. Third edition.
Graduate Texts in Mathematics, 20. Springer-Verlag,
New York, 1994.

\ii{[RRu]} J.-P.\ Rosay and W.\ Rudin:
Holomorphic maps from $\C^n$ to $\C^n$.
Trans.\ Amer.\ Math.\ Soc.\ {\bf 310}, 47--86 (1988)

\ii{[Sch]} M.\ Schneider:
On the number of equations needed to describe a variety.
Complex analysis of several variables (Madison, Wis., 1982),
163--180, Proc.\ Sympos.\ Pure Math., 41, Amer.\ Math.\ Soc.,
Providence, R.I., 1984.

\ii{[Sch\"ur]} J.\ Sch\"urmann:
Embeddings of Stein spaces into affine spaces of minimal dimension.
Math.\ Ann.\ {\bf 307}, 381--399 (1997).

\ii{[Whi]} G.\ W.\ Whitehead: Elements of homotopy theory.
Graduate Texts in Mathematics, 61. Springer, New York--Berlin, 1978.

%
%
%
%
\bigskip\medskip
\ni \sl Address: \rm 
IMFM, University of Ljubljana, Jadranska 19, 1000 Ljubljana, Slovenia

\bye